\documentclass[12pt]{amsart}
\input{amssymb.sty}
\title{Boundary Amenability of Relatively Hyperbolic Groups}
\author{Narutaka OZAWA}
\address{Department of Mathematics, UCLA, Los Angeles, 
CA 90095-1555\hspace*{\fill}\linebreak 
\indent Department of Mathematical Sciences,
University of Tokyo, Komaba, 153-8914}
\email{narutaka@ms.u-tokyo.ac.jp}
\date{April 1, 2005}
\thanks{The author was partly supported by JSPS}

\newtheorem{thm}{Theorem}
\newtheorem{prop}[thm]{Proposition}
\newtheorem{cor}[thm]{Corollary}
\newtheorem{lem}[thm]{Lemma}
\theoremstyle{definition}
\newtheorem{defn}[thm]{Definition}
\newcommand{\G}{\Gamma}
\newcommand{\dk}{\Delta K}
\newcommand{\bk}{\partial K}
\newcommand{\FF}{\mathcal{F}}
\newcommand{\e}{\varepsilon}
\newcommand{\supp}{\mathop{\mathrm{supp}}\nolimits}
\newcommand{\prob}{\mathop{\mathrm{Prob}}\nolimits}
\newcommand{\R}{{\mathbb R}}

\newcommand{\N}{{\mathbb N}} 
 
\newcommand{\ip}[1]{\langle#1\rangle} 
\begin{document}
\begin{abstract}
Let $K$ be a fine hyperbolic graph 
and $\G$ be a group acting on $K$ with finite quotient. 
We prove that $\G$ is exact provided that all vertex 
stabilizers are exact. 
In particular, a relatively hyperbolic group is exact if 
all its peripheral groups are exact. 
We prove this by showing that the group $\G$ acts amenably 
on a compact topological space. 
We include some applications to the theories 
of group von Neumann algebras and of measurable orbit 
equivalence relations.
\end{abstract}
\maketitle
\section{Introduction}
A discrete group is said to be \emph{exact} if it acts 
amenably on some compact topological space 
\cite{gh}\cite{exact}. 
Higson and Roe \cite{hr} proved that this property is equivalent 
to the property $A$ of Yu \cite{yu} and is a coarse equivalence 
invariant. 
They also proved that a group with finite asymptotic dimension 
is exact. 
It was proved by Adams \cite{adams} that the boundary action 
of a hyperbolic group on its ideal boundary is amenable and hence 
hyperbolic groups are exact. 
(See \cite{ar}\cite{kaimanovich} for simpler proofs.)
In this paper, we study the case of 
relatively hyperbolic groups. 
The notion of relative hyperbolicity was proposed 
by Gromov \cite{gromov} and has been developed 
in \cite{bowditch}\cite{ds}\cite{farb}\cite{osinmemoirs}, 
just to name a few.
In particular, Bowditch \cite{bowditch} introduced 
the notion of a fine hyperbolic graph $K$ and 
studied its compactification $\dk$. 
Bowditch characterized a relatively hyperbolic group 
as a group which admits a suitable action on 
a fine hyperbolic graph $K$. 
We will prove the following generalization of Adams's 
and Kaimanovich's theorems \cite{adams}\cite{kaimanovich}. 
\begin{thm}\label{thm}
Let $\G$ be a group acting on a countable fine hyperbolic graph $K$ 
with finite quotient. 
Let $Y$ be a compact space on which $\G$ acts. 
We assume that, for every vertex $x$ in $K$, 
the restricted action of the vertex stabilizer 
$\G^x$ on $Y$ is amenable. 
Then, the diagonal action of $\G$ on $\dk\times Y$ is amenable. 
\end{thm}

We include the Hausdorff property in the definition of compactness 
and assume that group actions on topological spaces are continuous. 
We note that if $\Lambda$ is an exact subgroup of $\G$, 
then the left multiplication action of $\Lambda$ on the 
Stone-\v Cech compactification $\beta\G$ of $\G$ is amenable. 
Since a tree is a fine hyperbolic graph, 
we obtain Dykema's theorem \cite{dykema} as a corollary. 
(See also \cite{tu}.)

\begin{cor}
Amalgamated free products of exact groups are exact. 
\end{cor}

Bowditch's characterization of relatively hyperbolic groups 
implies the following. 

\begin{cor}
Let $\G$ be a finitely generated group which is 
hyperbolic relative to a family of subgroups $\mathcal{G}$. 
Then, $\G$ is exact provided that 
each element of $\mathcal{G}$ is exact.
\end{cor}
Alternative proofs of this corollary are obtained by 
Dadarlat and Guentner \cite{dg} and Osin \cite{osin}. 
This result should be compared with Osin's result \cite{osin} 
on finiteness of asymptotic dimension.
Osin kindly pointed out that the results need not be true 
for weakly relatively hyperbolic groups; the boundedly 
generated universal finitely presented 
group constructed in \cite{osin} is 
not exact since there exists a finitely presented group 
that is not exact \cite{gromov2}.

Amenable action has several applications. 
Besides those follow from exactness 
(e.g., the Novikov conjecture \cite{yu}\cite{higson}), 
there are applications to the classification of orbit equivalence 
relations (this was the motivation of Adams \cite{adams}) 
as well as the classification of group von Neumann algebras. 
For instance, we show that a group which is hyperbolic relative to 
a family of amenable groups satisfies the property $(AO)$ defined 
in \cite{solid}
(and is in the class $\mathcal{C}$ defined in \cite{kurosh}).
Hence, the relevant results in \cite{solid}\cite{op}\cite{kurosh} 
are applicable to such $G$. 
Examples of such groups include the fundamental groups of complete non-compact 
finite-volume Riemannian manifolds with pinched negative sectional curvature 
(which are hyperbolic relative to nilpotent cusp subgroups) \cite{farb} and 
limit groups (which are hyperbolic relative to maximal non-cyclic 
abelian subgroups) \cite{dahmani2}\cite{alibegovic}. 
\section{Fine Hyperbolic Graphs}
In this section, we collect several facts on 
a fine hyperbolic graph $K$ and its compactification $\dk$. 
Most of them are taken from \cite{bowditch}. 

Let $K$ be a graph with vertex set $V(K)$ and edge set $E(K)$. 
(We allow no loops nor multiple edges). 
A \emph{path} of length $n$ connecting $x,y\in V$ is a 
sequence $x_0x_1\cdots x_n$ of vertices, with $x_0=x$ and 
$x_n=y$, and with each $x_i$ equal to or adjacent to $x_{i+1}$.
A \emph{circuit} is a closed path ($x_0=x_n$) such that 
all $x_0,\ldots,x_{n-1}$ are distinct. 
For a finite path $\alpha=x_0x_1\cdots x_n$, we set 
$\alpha_-=x_0$, $\alpha_+=x_n$ and $\alpha(k)=x_k$. 

We put a path metric $d$ on $V(K)$, where $d(x,y)$ is the length of 
shortest path in $K$ connecting $x$ to $y$. 
We assume the graph $K$ is connected so that 
$d(x,y)<\infty$ for every pair of vertices.
A finite or infinite path $\alpha$ is \emph{geodesic} 
if $d(\alpha(m),\alpha(n))=|m-n|$ for all $m,n$. 
For any three point $x,y,z\in V(K)$, we define the Gromov product by 
\[
\ip{x,y}_z=\frac{1}{2}(d(x,z)+d(y,z)-d(x,y)).
\]

\begin{defn}
We say the graph $K$ is \emph{hyperbolic} if there exists $\delta>0$ 
such that every geodesic triangle is $\delta$-thin; 
for any geodesic paths $\alpha$ and $\beta$ with $\alpha_-=\beta_-=z$ 
and any $k\le\ip{\alpha_+,\beta_+}_z$, we have
$d(\alpha(k),\beta(k))<\delta$. 

We say the graph $K$ is \emph{fine} if 
for any $n$ and any edge $e\in E(K)$, 
the set $C(e,n)$ of circuits of length at most $n$ 
containing the edge $e$ is finite. 
\end{defn}

From now on, we assume that $K$ is a countable hyperbolic graph 
which is \emph{uniformly fine}, i.e., for any $n$, 
we have $\sup_{e\in E(K)}|C(e,n)|<\infty$. 
We note that any fine graph which admits a group action 
with finite quotient is uniformly fine. 
Two infinite geodesic paths $\alpha$ and $\beta$ are \emph{equivalent} if 
their Hausdorff distance is finite. 
The Gromov boundary $\bk$ of $K$ is defined as the set of 
all equivalence classes of infinite geodesic paths in $K$. 
We write $\dk=V(K)\cup\bk$. 
For an infinite geodesic path $\alpha=x_0x_1\cdots$, 
we denote by $\alpha_+$ the boundary point that is represented by $\alpha$. 
For a biinfinite geodesic path $\alpha=\cdots x_{-1}x_0x_1\cdots$, 
we likewise denote by $\alpha_-$ the boundary point that is represented by 
the geodesic path $x_0x_{-1}x_{-2}\cdots$. 
In any case, we say $\alpha$ connects $\alpha_-$ to $\alpha_+$. 
For every $x,y\in\dk$, we denote by $\FF(x,y)$ 
the set of all geodesic paths which connects $x$ to $y$.
\begin{lem}
For any $x,y\in\dk$, we have $\FF(x,y)\neq\emptyset$. 
\end{lem}
For every $x\in\dk$ and a finite subset $A\subset V(K)$, we define 
\[
M(x,A)=\{ z\in\dk : \mbox{ any $\alpha\in\FF(x,z)$ 
does not intersects with $A\setminus\{x\}$}\}
\]
The following is Bowditch's theorem (Section 8 in \cite{bowditch}).
\begin{thm}\label{thm:bowditch}
The family $\{M(x,A)\}_A$ defines a neighborhood base for $x\in\dk$. 
With this topology, $\dk$ is a compact topological space, 
in which $V(K)$ is dense. 
Every graph automorphism on $K$ uniquely extends to a homeomorphism on $\dk$.
\end{thm}
We note that $\bk$ is a Borel subset in $\dk$ since it is $G_\delta$. 
We claim that every $M(x,A)$ is open. 
For simplicity, we just prove that 
for every $x,y\in V(K)$, the set 
\[
T(x,y)=M(x,\{y\})^c
=\{ z\in\dk : \exists\alpha\in\FF(x,z)\mbox{ such that }y\in\alpha\}
\]
is closed. 
Let $(z_n)_n$ be a sequence in $T(x,y)$ 
which converges to $z$ in $\dk$. 
Let $\alpha_n$ be a geodesic path connecting 
$x$ to $z_n$ which passes on $y$. 
Let $l$ be the largest integer (possibly $\infty$) such that 
there exists $w\in V(K)$ with $\{n : \alpha_n(l)=w\}$ is infinite. 
If $l$ is finite, then $w$ is a limit point of the sequence $(z_n)_n$. 
(To see this, use the equivalent topology defined by $\{M'(w,B)\}_B$ 
in Section 8 in \cite{bowditch}.) 
This implies that $z=w\in T(x,y)$. 
Now, suppose that $l$ is infinite. 
Since the set of geodesic paths connecting any two points 
is locally finite, we may pass to a subsequence and assume that 
$\alpha_n(1)=\cdots=\alpha_n(n)$ for every $n$. 
It follows that the geodesic path $\alpha$, given by $\alpha(n)=\alpha_n(n)$, 
connects $x$ to $z$ and hence $z\in T(x,y)$. 

For every $x\in V(K)$, $z\in\bk$ and $l,k\in\N$, we set 
\[
S(x,z,l,k)=\{ \alpha(l)\in V(K) : \alpha\in\FF(x',z)\mbox{ with }
d(x',x)\le k\}.
\]
We observe that $S(x,z,l,k)$ is finite. 
For a finite subset $S\subset V(K)$, we write $\xi_S=|S|^{-1}\chi_S$ 
for the normalized characteristic function on $S$. 
We note that for every $x,y\in V(K)$ and $l,k\in\N$, the set 
\[
\{ z\in\bk : y\in S(x,z,l,k)\}=
\bk\cap\bigcup_{x'}\{T(x',y) : d(x',x)\le k,\ d(x',y)=l\}
\]
is Borel in $\dk$. 
It follows that, for every $x,y\in V(K)$ and $l,k\in\N$, the functions 
\[
\bk\ni z\mapsto |S(x,z,l,k)|\in\N
\quad\mbox{ and }\quad
\bk\ni z\mapsto\xi_{S(x,z,l,k)}(y)\in\R
\]
are Borel. 
The following lemma is well-known for 
uniformly locally finite hyperbolic graph. 
It also follows from Lemma~16 of \cite{osin} 
in special cases.  

\begin{lem}\label{lem:lin}
There exists a constant $C=C(K)>0$ such that 
$|S(x,z,l,k)|<Ck$ for all $x\in V(K)$, $z\in\bk$ 
and $l,k\in\N$ with $l>k+\delta$. 
\end{lem}
\begin{proof}
Let $\alpha$ be any fixed geodesic path connecting $x$ to $z$. 
It suffices to show 
\[
S(x,z,l,k)\subset \bigcup \{C(e,6\delta) : e \mbox{ an edge in }\alpha([l-k-\delta,l+k+\delta])\},
\]
where $C(e,n)$ is the set of circuits of length at most $n$ 
containing the edge $e$. 
Choose $y\in S(x,z,l,k)$. 
There exists a geodesic path $\beta$ connecting $x'$ to $z$ 
such that $d(x',x)\le k$ and $y=\beta(l)$. 
Since any geodesic triangle is $\delta$-thin, 
there exists a number $m\in[l-k,l+k]$ such that 
$d(\alpha(m\pm\delta),\beta(l\pm\delta))<\delta$. 
Let $\gamma^{\pm}$ be any geodesic paths connecting 
$\alpha(m\pm\delta)$ to $\beta(l\pm\delta)$. 
Since $d(y,\gamma^+\cup\gamma^-)\geq1$ and $d(\gamma^+,\gamma^-)\geq2$, 
the path $\gamma^-\cup\beta([l-\delta,l+\delta])\cup\gamma^+$ 
do not contain a circuit which includes $y$. 
Thus, the closed path 
$\gamma^-\cup\beta([l-\delta,l+\delta])\cup\gamma^+\cup\alpha([m-\delta,m+\delta])$ 
of length at most $6\delta$ contains 
a circuit which includes both $y$ and 
an edge in $\alpha([m-\delta,m+\delta])$. 
\end{proof}
For every $x\in V(K)$, $z\in\bk$ and $n$, we set 
\[
\zeta^n_{x,z}=\frac{1}{n}\sum_{k=n+1}^{2n}\xi_{S(x,z,4n,k)}\in\ell_1(V(K)).
\]
We note that the map $\bk\ni z\mapsto\zeta^n_{x,z}(y)\in\R$ is 
Borel for every $x,y\in V(K)$.
The following lemma is borrowed from \cite{kaimanovich}, but we put a proof for completeness.
\begin{lem}\label{lem:ave}
For any $d>0$, we have 
\[
\lim_{n\to\infty}\sup_{z\in\bk}\sup_{d(x,x')\le d}
\|\zeta^n_{x,z}-\zeta^n_{x',z}\|_1 =0.
\]
\end{lem}
\begin{proof}
Let $z\in\bk$ and $x,x'\in V(K)$ be given. 
Let $d=d(x,x')$ and $n\geq d+\delta$. 
For simplicity, we write $S_k=S(x,z,4n,k)$ and $S'_k=S(x',z,4n,k)$. 
Then, we have
$S_k\cup S'_k\subset S_{k+d}$ and $S_k\cap S'_k\supset S_{k-d}$
for every $n< k\le 2n$.
It follows that 
\[
\|\xi_{S_k}-\xi_{S'_k}\|_1
=2\left(1-\frac{|S_k\cap S'_k|}{\max\{|S_k|,|S'_k|\}}\right)
\le 2\left(1-\frac{|S_{k-d}|}{|S_{k+d}|}\right)
\]
for $n<k\le 2n$. 
Since $|S_k|\le Ck$ by Lemma~\ref{lem:lin}, we have 
\[
\frac{1}{n}\sum_{k=n+1}^{2n}
 \frac{|S_{k-d}|}{|S_{k+d}|}\geq
\left(\prod_{k=n+1}^{2n}
 \frac{|S_{k-d}|}{|S_{k+d}|}\right)^{1/n}
=\left(
 \frac{\prod_{k=n+1-d}^{n+d}|S_k|}{\prod_{k=2n+1-d}^{2n+d}|S_k|}
 \right)^{1/n}\geq 
(3Cn)^{-2d/n}.
\]
It follows that 
\[
\|\zeta^n_{x,z}-\zeta^n_{x',z}\|_1
\le\frac{1}{n}\sum_{k=n+1}^{2n}\|\xi_{S_k}-\xi_{S'_k}\|_1 
\le 2(1-(3Cn)^{-2d/n}).
\]
Since $(3Cn)^{-2d/n}\to 1$ as $n\to\infty$, 
this proves the lemma.
\end{proof}
\begin{lem}\label{lem:infty}
Let $(x_n)_n$ and $(y_n)_n$ be sequences in $V(K)$ 
which converge to respectively $x$ and $y$ in $\dk$. 
Suppose that $\sup_n d(x_n,y_n)<\infty$. 
Then, we have either that $x=y$, or that 
for every sufficiently large $n$, 
all geodesic paths connecting $x_n$ to $y_n$ 
intersect with both $x$ and $y$. 
\end{lem}
\begin{proof}
We observe first that if $x\in\bk$, then $y\in\bk$ and $x=y$. 
Thus, we assume $x,y\in V(K)$. 
Let $C=\sup_n d(x_n,y_n)$. 
We claim that if $\limsup d(x,x_n)=\infty$, then $x=y$. 
Indeed, passing to a subsequence, we may assume 
that $d(x,x_n)>10C$ for all $n$. 
For any finite subset $A\subset V(K)\setminus\{x\}$, 
there exists $\alpha_n\in\FF(x,x_n)$ which does not intersect with $A$. 
Composing $\alpha_n$ with a geodesic path connecting $x_n$ to $y_n$ 
and then discarding redundant vertices, we obtain a 
$10C$-quasigeodesic path connecting $x$ to $y_n$ 
which does not intersects with $A$. 
Hence, the sequence $(y_n)_n$ converges to $x$. 
(To see this, use the equivalent topology defined by $\{M'_{10C}(x,A)\}_A$ 
in Section 8 in \cite{bowditch}). 
This proves the claim. 
Finally, we suppose that $C'=\sup d(x,x_n)<\infty$. 
Let $\alpha_n\in\FF(x,x_n)$, $\beta_n\in\FF(y,y_n)$ and 
$\gamma\in FF(x,y)$. 
Since $x_n\to x$, we have $\alpha_n\cap\gamma=\{x\}$ for sufficiently large $n$. 
Thus any geodesic path connecting $x_n$ to $y_n$, which does not 
intersects $x$, gives rise to a circuit of length at most $2(C+C'+d(x,y))$ 
containing $x$ and edges in $\alpha_n$ and $\gamma$. 
Since $x_n\to x$, this is only possible for finitely many $n$. 
The same thing for $y$ and we are done. 
\end{proof}
\section{Amenable Actions}
For every locally compact space $\Omega$, 
we denote by $\prob(\Omega)$ the space of 
all regular Borel probability measure on $\Omega$, 
equipped with the weak$^*$-topology. 
In particular, if $\G$ is a discrete group, then 
\[
\prob(\G)=\{ \mu \in\ell_1(\G) : 
\mu\geq0\mbox{ and }\sum_{s\in\G}\mu(s)=1\}
\]
and $\prob(\G)$ is equipped with the pointwise convergence 
topology. 
We note that the pointwise convergence topology coincides 
with the norm topology on $\prob(\G)$.
The group $\G$ acts on $\prob(\G)$ from the left; 
$(s\cdot\mu)(t)=\mu(s^{-1}t)$ for every $s,t\in\G$ and $\mu\in\prob(\G)$. 
Let $\Omega$ be a locally compact space. 
By an \emph{action} of $\G$ on $\Omega$, we mean a group homomorphism 
from $\G$ into the group of homeomorphisms on $\Omega$. 
\begin{defn}\label{defn:amenable}
Let $\G$ be a discrete group acting 
on a compact topological space $\Omega$. 
We say the action of $\G$ on $\Omega$ is \emph{amenable} 
if for every finite subset $E\subset\G$ and $\e>0$, 
there exists a continuous map  
\[
\mu\colon\Omega\ni\omega\mapsto\mu_\omega\in\prob(\G)
\]
such that 
\[
\max_{s\in E}\sup_{\omega\in\Omega}
\|s\cdot\mu_\omega-\mu_{s\omega}\|\le\e.
\]
We say a group $\G$ is \emph{exact} if it acts amenably on 
some compact topological space. 
\end{defn}
A group $\G$ is amenable iff the trivial action 
on a singleton set is amenable. 
If a group $\G$ is exact, then it acts amenably on the Stone-\v Cech 
compactification $\beta\G$. 
We will use the following application of 
the Hahn-Banach theorem in the spirit of Anantharaman-Delaroche \cite{delaroche}\cite{ar}. 
If $\G$ acts on $K$, we denote the stabilizer subgroup 
of $a\in K$ by $\G^a=\{ s\in \G : sa=a\}$. 

\begin{prop}\label{prop:borel}
Let $\G$ be a countable group acting on $X$, $Y$ and $K$, 
where $X$ and $Y$ are compact and $K$ is countable discrete. 
Assume that for any finite subset $E\subset\G$ and $\e>0$ 
there exists a Borel map 
$\zeta\colon X\to\prob(K)$ 
(i.e., the function $X\ni x\mapsto\zeta_x(a)\in\R$ is 
Borel for every $a\in K$) such that 
\[
\max_{s\in E}\sup_{x\in X}\| s\zeta_x-\zeta_{sx}\|<\e.
\]
Assume moreover that the restricted action of 
the vertex stabilizer $\G^a$ on $Y$ is amenable for every $a\in K$. 
Then, the diagonal action of $\G$ on $X\times Y$ is amenable. 
\end{prop}
\begin{proof}
We first claim that we may take $\zeta$ in the statement 
to be continuous rather than Borel. 
Fix a finite symmetric subset $E\subset\G$. 
For every continuous map $\zeta\colon X\to\prob(K)$, 
we define $f_\zeta\in C(X)$ by 
\[
f_\zeta(x)=\sum_{s\in E}\|s\zeta_x-\zeta_{sx}\|
=\sum_{s\in E}\sum_{a\in K}|\zeta_x(s^{-1}a)-\zeta_{sx}(a)|.
\]
Since $f_{\sum_k\alpha_k\zeta_k}\le\sum_k\alpha_kf_{\zeta_k}$ 
for every $\alpha_k\geq 0$ with $\sum_k\alpha_k=1$, 
if $0$ is in the weak closure of $\{ f_\zeta : \zeta\}$ in $C(X)$, 
then $0$ is in the norm closure of $\{ f_\zeta :\zeta\}$ by 
the Hahn-Banach separation theorem. 
We note that the dual of $C(X)$ 
is the space of finite regular Borel measures on $X$ 
by the Riesz representation theorem. 
Let $\e>0$ and $m\in\prob(X)$ be given.
By assumption, 
there exists a Borel map $\eta\colon X\to\prob(K)$ 
such that 
$\sup_{x\in X}\| s\eta_x-\eta_{sx}\|<\e/|E|$.
By the countable additivity of the measure, 
there exists a finite subset $F\subset K$ such that 
$\int_X \sum_{a\in F}\eta_x(a)\,dm(x)>1-\e/|E|$. 
We approximate, for each $a\in F$, the Borel function $x\mapsto\eta_x(a)$ 
by a continuous function and obtain a continuous map 
$\zeta\colon X\to\prob(K)$ such that 
$\supp\zeta_x\subset F$ for all $x\in X$ and 
\[
\int_X\|\zeta_x-\eta_x\|\,dm(x)=\int_X\sum_{a\in K}|\zeta_x(a)-\eta_x(a)|\,dm(x)<2\e/|E|.
\]
It follows that 
\[
\int_X f_\zeta(x)\,dm(x)
<\int_X \sum_{s\in E}\|s\eta_x-\eta_{sx}\|\,dm(x)+4\e<5\e.
\]
Thus, we proved our claim. 

Now, let a finite subset $E\subset\G$ and $\e>0$ be given. 
By the previous result, there exists a continuous $\zeta$ 
such that $\sup_{x\in X}\|s\zeta_x-\zeta_{sx}\|<\e$ for every $s\in E$. 
We may assume that there exists a finite subset $F\subset\G$ such that 
$\supp\zeta_x\subset F$ for all $x\in X$. 
We fix a $\G$-fundamental domain $V\subset K$ with 
a projection $v\colon K\to V$ and 
a cross section $\sigma\colon K\to\G$, i.e., 
$K$ decomposes into the disjoint union $\bigsqcup_{v\in V}\G v$ 
and $a=\sigma(a)v(a)$ for every $a\in K$. 
We note that $\sigma(sa)^{-1}s\sigma(a)\in\G^{v(a)}$ 
for every $s\in\G$ and $a\in K$. 
For each $v\in V$, we set
\[
E^v=\{ \sigma(sa)^{-1}s\sigma(a) : a\in F\cap\G v\mbox{ and }s\in E\}
\subset\G^v.
\]
Since the $\G^v$ action on $Y$ is amenable and $E^v$ is finite, there exists 
a continuous map $\nu^v\colon Y\to\prob(\G)$ such that 
\[
\max_{s\in E^v}\sup_{y\in Y}\|s\nu^v_y-\nu^v_{sy}\|<\e.
\]
Now, we define $\mu\colon X\times Y\to\prob(\G)$ by
\[
\mu_{x,y}=\sum_{a\in K}\zeta_x(a)\,\sigma(a)\nu^{v(a)}_{\sigma(a)^{-1}y}.
\]
The map $\mu$ is clearly continuous. 
Moreover, we have
\begin{alignat*}{2}
s\mu_{x,y} &=  \sum_{a\in K}\zeta_x(a)\, s\sigma(a)\nu^{v(a)}_{\sigma(a)^{-1}y} &
&= \sum_{a\in K}\zeta_x(a)\, \sigma(sa)\big(\sigma(sa)^{-1}s\sigma(a)\nu^{v(a)}_{\sigma(a)^{-1}y}\big)\\
& & &\approx_\e \sum_{a\in K}\zeta_x(a)\, \sigma(sa)\nu^{v(sa)}_{\sigma(sa)^{-1}sy}\\
& & &\approx_\e \sum_{a\in K}\zeta_{sx}(sa)\, \sigma(sa)\nu^{v(sa)}_{\sigma(sa)^{-1}sy}
=\mu_{sx,sy}
\end{alignat*}
for every $s\in E$ and $(x,y)\in X\times Y$. 
\end{proof}
We prove Theorem~\ref{thm}. 
\begin{proof}
We may assume that $\G$ is countable 
and $K$ is uniformly fine. 
Let a finite subset $E\subset\G$ and $\e>0$ be given. 
Fix an origin $o\in V(K)$.
By Lemma~\ref{lem:ave}, there exists $n$ such that 
$\zeta\colon\bk\ni z\mapsto\zeta_z=\zeta^n_{o,z}\in\prob(K)$ satisfies 
\[
\max_{s\in E}\sup_{z\in\bk}\|s\zeta_z-\zeta_{sz}\|
=\max_{s\in E}\sup_{z\in\bk}\|\zeta^n_{so,sz}-\zeta^n_{o,sz}\|<\e.
\]
For $x\in V(K)$, we set $\zeta_x=\delta_x\in\prob(K)$. 
Then, the map $\zeta\colon\dk\to\prob(K)$ is Borel (see the remark 
preceding Lemma~\ref{lem:ave}) 
and satisfies the assumptions given in Proposition~\ref{prop:borel}. 
Thus the action on $\dk\times Y$ is amenable.
\end{proof}
\section{Applications}
Recall that a \emph{compactification} of a discrete group $\G$ 
is a compact topological space $\bar{\G}$ which contains 
$\G$ as a discrete open dense subset. 
We assume that the left translation action of $\G$ on $\G$ 
extends to a continuous action of $\G$ on $\bar{\G}$. 
We have well-known one-to-one correspondence between 
such compactification $\bar{\G}$ and a C$^*$-subalgebra 
$c_0(\G)\subset C(\bar{\G})\subset\ell_\infty(\G)$ 
which is invariant under left translation. 
We say the compactification $\bar{\G}$
is \textit{small at infinity} if 
the right translation action of $\G$ on $\G$ extends 
to a continuous action on $\bar{\G}$ and if its restriction 
to $\partial\G=\bar{\G}\setminus\G$ is trivial. 
In other words, $f^t-f\in c_0(\G)$ 
for all $f\in C(\bar{\G})\subset\ell_\infty(\G)$ and $t\in\G$, 
where $f^t(s)=f(st^{-1})$ for $s\in\G$. 

\begin{prop}\label{prop}
Let $\G$ be a group which is hyperbolic relative to 
a family $\mathcal{G}$ of amenable subgroups. 
Then, $\G$ acts amenably on some compactification $\bar{\G}$ 
which is small at infinity. 
\end{prop}
\begin{proof}
By \cite{bowditch}, the group $\G$ admits a finite quotient 
action on a fine hyperbolic graph $K$ in a way that 
every infinite vertex stabilizer is in $\mathcal{G}$ 
and that $\G^x\cap\G^y$ is finite for any $x,y\in V(K)$ with $x\neq y$.  
By Theorem~\ref{thm}, the action of $\G$ on $\dk$ is amenable. 
We fix an origin $o$ and consider the map 
$\G\ni s\mapsto so\in \dk$. 
This map induces a compactification $\bar{\G}$ 
such that $\bar{\G}\to\dk$ is continuous. 
We note that $\G$ acts on $\bar{\G}$ amenably.   
To prove $\bar{\G}$ is small at infinity, it suffices to show 
that for any $t\in\G$ and $f\in C(\dk)$, we have 
$f_o-f_o^t\in c_0(\G)$, where $f_o\in\ell_\infty(\G)$ 
is defined by $f_o(s)=f(so)$ for $s\in\G$. 
We fix $\e>0$ and show that 
$A=\{ s\in\G : |f(so)-f(st^{-1}o)|\geq\e\}$ is finite. 
Let $(s_n)_n\subset A$ be any sequence such that 
all $s_n$ are distinct. 
Since $\dk$ is compact (and first countable), 
we may assume that $s_no\to x$ and $s_nt^{-1}o\to y$. 
We note that $x\neq y$ since $|f(x)-f(y)|\geq\e$. 
Let $\alpha$ be a geodesic path connecting $o$ to $t^{-1}o$. 
Then, by Lemma~\ref{lem:infty}, for all sufficiently large $n$, 
the geodesic paths $s_n\alpha$ intersect with both $x$ and $y$. 
Since $\G^x\cap\G^y$ is finite, $A$ is finite. 
\end{proof}
It follows from (the proof of) Lemma 5.2 in \cite{hg} that 
the group $\G$ satisfying the assumptions of Proposition~\ref{prop} 
is in the class $\mathcal{C}$ defined in \cite{kurosh} 
(i.e., the left and right translation action of $\G\times\G$ 
on the Stone-\v Cech remainder $\partial^\beta\G$ is amenable). 
Hence, the corresponding results in \cite{solid}\cite{op}\cite{kurosh} 
are applicable to such $\G$.

\end{document}